\documentclass[a4paper,12pt]{article} 
\newdimen\paperhight
\usepackage{amssymb}
\usepackage{amsfonts}
\usepackage[usenames]{color}

\newcommand{\dsp}{\displaystyle }

\newcommand{\hf}{\frac{1}{2}}
\newcommand{\st}{\frac{1}{16}}
\newcommand{\tth}{\frac{2}{3}}
\newcommand{\ff}{\frac{4}{5}}
\newcommand{\ssv}{\frac{6}{7}}

\newcommand{\pr}{\par \vspace{3mm}\noindent [{\bf Proof}] \qquad}
\newcommand{\prend}{\hfill \qed \par \vspace{3mm}}

\newcommand{\qed}{\quad \hbox{\rule[-2pt]{3pt}{6pt}} \par \vspace{3mm}}

\newcommand{\1}{{\bf 1}} 
 
\newcommand{\C}{\mathbb C} 
\newcommand{\Z}{\mathbb Z}

\newcommand{\M}{\mathbb M} 
\newcommand{\N}{\mathbb N} 
\newcommand{\R}{\mathbb R} 
\newcommand{\CG}{{\cal G}}

\newcommand{\al}{\alpha}
\newcommand{\be}{\beta}

\newcommand{\ep}{\epsilon}
\newcommand{\ga}{\gamma}
\newcommand{\la}{\lambda}

\newtheorem{thm}{Theorem}[section]
\newtheorem{prn}[thm]{Proposition}

\newtheorem{lmm}[thm]{Lemma}

\newtheorem{cnj}{Conjecture}

\begin{document}
\title{VOAs generated by two conformal vectors whose 
$\tau$-involutions generate $S_3$.}
\author{Masahiko Miyamoto 
\footnote{Supported by the Grants-in-Aids for Scientific Research, 
No. 13440002 and No. 12874001, The Ministry of Education, 
Science and Culture, Japan.}} 
\date{\begin{tabular}{c}
Institute of Mathematics \cr
University of Tsukuba \cr
Tsukuba 305, Japan \cr
\end{tabular}}
\maketitle

\begin{center}
{\large Dedicated to Professor Koichiro Harada for his 60th 
birthday.} 
\end{center}

\begin{abstract}
We determined the inner products of two conformal vectors with 
central charge $\hf$ whose $\tau$-involutions generates $S_3$ if 
none of $\tau$-involutions are trivial. 
We also see that a subVA generated by such conformal vectors 
is a VOA with central charge $\hf+\frac{20}{21}$ or 
$\frac{4}{5}+\frac{6}{7}$ and has a 
Griess algebra of dimension three or four, respectively.
\end{abstract}

\section{Introduction}
A vertex operator algebra (shortly VOA) is a mathematical object 
for a 2-dimensional conformal field theory, but it comes from the 
Moonshine conjecture to explain mysterious properties of 
Monster simple group $\M$ \cite{CN}. As an answer to this conjecture, 
Frenkel, Lepowsky and Meurman construct the moonshine VOA 
$V^{\natural}=\oplus_{i=0}^{\infty}V^{\natural}_i$ 
whose full automorphism group is the Monster simple group $\M$, 
\cite{FLM}. 
Its weight two subspace $V^{\natural}_2$ coincides with a 
commutative (non-associative) algebra (called the monstrous Griess 
algebra) 
of dimension 196884 
constructed by Griess in order to construct the Monster simple 
group \cite{Gr}. 
This algebra is studied from a group theoretic point of view. One of 
the important results is that each $2A$-involution $\theta$ defines 
a unique idempotent $e_{\theta}$ (called an axis) of the 
monstrous Griess algebra such that the inner product 
$\langle e_{\theta}, e_{\phi}\rangle$ is uniquely determined 
by the conjugacy classes of 
$\theta\phi$, see \cite{Co}. 

From the view point of vertex operator algebras, 
the author showed that an involutive automorphism $\tau_e$ comes 
arise from 
a rational conformal vector with central charge $\hf$ 
in a Griess algebra, namely, 
if $e$ generates a rational VOA $L(\hf,0)$ called an Ising model. 
then one can define  
an involutive automorphism $\tau_e$ of $V$ by 
$$  \tau_e: \quad \left\{
\begin{array}{rl}
1  &\mbox{ on  } W_0\oplus W_{{1\over 2}} \cr
-1  &\mbox{ on  } W_{{1\over 16}},  
\end{array} \right.  $$
where $W_h$ denotes the sum of all irreducible ${\rm VA}(e)$-modules 
isomorphic to $L(\hf,h)$ and ${\rm VA}(e)$ is a subVOA generated 
by $e$.  

In the monstrous Griess algebra, it is easy to check that 
such a conformal vector $e$ (with central charge $\hf$) 
is corresponding to an axis and $\tau_e$ is a $2A$-involution. 
Even in a (general) Griess algebra, the inner products are 
very important and the author showed that if two involutions 
$\tau_e$ and $\tau_f$ commute, that is, $\tau_e\tau_f$ is of 
order two,  
then the inner product $\langle e,f\rangle $ is $0$ or $\frac{1}{32}$. 
In the moonshine VOA $V^{\natural}$, 
$\tau_e\tau_f$ is a $2A$-involution if and only if 
$\langle e,f\rangle =\frac{1}{32}$ 
and a $2B$-involution if and only if $\langle e,f\rangle=0$. 
Conversely, if $\langle e,f\rangle=0$, then $\tau_e$ and $\tau_f$ 
commute. 
One of the mysteries about the Monster simple group is 
a property of $2A$-involutions. 
The $2A$-involutions satisfy several interesting properties.  
For example, they satisfy a $6$-transposition property, that is, 
$|\tau_e\tau_f|\leq 6$ for any $2A$-involutions $\tau_e$ and $\tau_f$. 
There is also a mysterious relation with  
$E_8$ Dynkin diagram, which was observed by McKay. Namely, the 
conjugacy classes of $\tau_e\tau_f$ for distinct 2A-involutions 
$\tau_e$ and 
$\tau_f$ are \\
$$ 2A, \ 3A, \ 4A, \ 5A, \ 6A, \ 3C, \ 4B \ \mbox{  and  }2B $$
and these numbers are multiplicities of roots in a maximal root in 
$E_8$-root system. 
Another interesting topic is 
a $Y_{5,5,5}$-diagram or the set of $26$ involutions whose graph 
coincides with the incident graph 
given by $13$ lines and $13$ points of the projective plane of 
order three and the Bimonster $\M\wr \Z_2$ contains such a set as 
generators,  
\cite{At}. 
If we restrict them to the Monster simple group $\M$, which  
contains $21$ $2A$-involutions whose graph is 
the incident graph of $12$ lines and $9$ points in the Affine 
plane of 
order three, where a vertex is a $2A$-involution $\tau_e$ and 
an edge $\tau_e -\tau_f$ means 
$|\tau_e\tau_f|=3$ and no edge between $\tau_e$ and $\tau_f$ 
implies that $\tau_e\tau_f$ is of order two, see \cite{Mi1}. 
Since the author has already studied the case without edge in 
\cite{Mi2}, 
it becomes very important to treat the case with an edge, that is, 
the relation between two conformal vectors $e$ and $f$ 
such that $\tau_e\tau_f$ is of order three. \\

In this paper, we will study such a pair $(e, f)$ in a (general) 
Griess algebra with a positive definite invariant bilinear form 
$\langle\cdot,\cdot\rangle$. 
We also add one more assumption that none of $\tau$-involutions 
are trivial. 
We should note that if $\tau_e=1$, then we can define an automorphism 
$\sigma_e$ of different type by 
$$  \tau_e: \quad \left\{
\begin{array}{rl}
1  &\mbox{ on  } W_0 \cr
-1  &\mbox{ on  } W_{{1\over 2}}  
\end{array} \right.  $$
and 
$<\sigma_e| e\in {\rm Con}(\hf), \tau_e=1>$ is a normal 
3-transposition subgroup of the full 
automorphism group of $V$ and such 3-transposition groups 
are studied by \cite{F}, \cite{GH} and \cite{KMi}, where 
${\rm Con}(\hf)$ 
is the set of all conformal vector with central charge $\hf$. 
Moreover, if $\sigma_e$ is also trivial, then $V$ is a tensor 
product of a sub VOA ${\rm VA}(e)$ and some subVOA $W$. 
Since we are interested in the property of group generated by 
$\tau$-involutions, we get rid of these cases and assume that 
$\tau_e\not=1$ for any conformal vector $e$. \\

In the Monster group, if a 
product of two $2A$-involutions has order three, then it is a 
$3A$-triality 
or a $3C$-triality. For such conformal vectors in the monstrous Griess 
algebra $V^{\natural}_2$, the inner products are 
$\frac{13}{2^{10}}$ and 
$\frac{1}{2^8}$ corresponding to $3A$-triality or $3C$-triality, 
respectively, see \cite{Co}. 
We will show that this result is true for general 
VOAs if none of $\tau$-involutions are trivial, 
that is, there are only two possibilities of inner products
$\langle e,f\rangle $. Moreover, we will study a 
subVA ${\rm VA}(e,f)$ generated by $e$ and $f$. 
For example, we will show that ${\rm VA}(e,f)$ is 
a vertex operator algebra. 
We note that the notation ${\rm VA}(e,f)$ denotes a vertex subalgebra 
generated by $e$ and $f$ and we don't assume the existence of 
Virasoro element. 
Since each edge in the graph of 21 involutions 
is corresponding to a $3A$-triality, 
we focus our attention on the case where 
the inner product is $\frac{13}{2^{10}}$. 
In this case, ${\rm VA}(e,f)$ is a VOA with 
central charge $\frac{58}{35}$ and $\dim ({\rm VA}(e,f))_2=4$. 
Moreover, ${\rm VA}(e,f)$ contains 
$L(\hf,0)\otimes L(\frac{81}{70},0)$ and also 
$L(\ff,0)\otimes L(\ssv,0)$, where $L(c,0)$ denotes a simple 
Virasoro 
VOA with central charge $c$ and $L(c,h)$ denotes its irreducible 
module 
with highest weight $h$. 
In particular, a Griess algebra of ${\rm VA}(e,f)$ coincides with 
the subalgebra generated of $e$ and $f$ as a subalgebra of Griess 
algebra. 
We will call this a VOA of $\tau$-involution type $A_2$. 
We may view $L(\hf,0)$ as a VOA of 
involution type $A_1$.

\noindent 
{\bf  Acknowledgment}  \\ 
The author wishes to thank R. L. Griess for his helpful advices. 

\section{Setting and products in Griess algebra}
This paper is a continuation of \cite{Mi2} and we will adopt 
the notation from it. 
Since our interest is a finite automorphism group, we will treat a 
simple 
VOA $(V,Y,\1,\omega)$ over the real number field ${\R}$ 
and $\C V$ denotes its complexification $\C\otimes_{\R}V$. 
For $v\in V$, 
$Y(v,z)=\sum_{n\in \Z}v(n)z^{-n-1}$ denotes the vertex operator 
of $v$. \\

\noindent
Assume the following conditions: \\
(1) $V=\oplus_{n=0}^{\infty}V_n$, \ $V_0=\R \1$  and \ $V_1=0$. \\

A VOA of this type is called OZ(one zero)-type. 
Since $\dim V_0=1$ and $V_1=0$, 
there is a unique invariant bilinear form 
$\langle \cdot,\cdot\rangle$ 
on $V$ satisfying $\langle \1,\1\rangle=1$. We also assume \\

\noindent
(2) $\langle\cdot,\cdot\rangle$ is definite on $V_n$ for 
each $n$. \\

In particular, $V_2$ has a positive definite bilinear 
form $\langle\cdot,\cdot\rangle$ satisfying 
$\langle a(1)b,c\rangle=\langle b,a(1)c\rangle$ 
for $a,b,c\in V_2$, 
where $\langle a,b\rangle$ is given by 
$\langle a,b\rangle\1=a(3)b\in \R \1$.
A VOA of this type is called a moonshine type. 
When we consider a VOA $\C V$, we define the inner product on 
$\C\otimes_{\R}V_2$ by 
$\langle a,b\rangle\1=a_3\overline{b}$, where $\overline{b}$ is a 
complex conjugate of $b$. 

Under the above assumptions, 
$(V_2, \cdot(1)\cdot)$ becomes a commutative (non-associative) 
algebra 
called a Griess algebra with a positive definite invariant 
bilinear form. To simplify the notation, 
$ef$ denotes $e(1)f$ for $e,f\in V_2$. \\

As we mentioned in the introduction, we add the following assumption. 

\noindent
(3) For any conformal vector $e$ with central charge $\hf$, 
$\tau_e$ is not trivial. \\

In this case, there is a one-to-one correspondence between 
$\tau$-involutions and conformal vectors with central charge $\hf$, 
see \cite{Mi2}.

Since $V$ has a definite invariant bilinear form, if $V$ contains 
a Virasoro VOA $L(c,0)$, $L(c,0)$ has a unitary highest 
weight representation in $V$. The work in \cite{FQS} and \cite{GKO} 
gives a 
complete classification 
of unitary highest weight representations of the Virasoro algebra. 
They proved that the highest weight representation $L(c,h)$ 
is unitary if and only if 
either $(c,h)$ satisfies $c\geq 1$ and $h\geq 0$, or else $(c,h)$ 
is among the following list:
$$ \begin{array}{c}
\dsp{  c=c_m=1-\frac{6}{(m+2)(m+3)} \quad (m=0,1,2,...),} \cr 
   \cr
\dsp{h=h^m_{r,s}=\frac{[(m+3)r-(m+2)s]^2-1}{4(m+2)(m+3)}} \quad 
(r,s\in {\N}, 1\leq s\leq r\leq m+1).
\end{array}  \eqno{(2.1)}$$
The unitary representations $L(c_m,h^m_{r,s})$ for $(c_m, h^m_{r,s})$ in 
the discrete series as above are called the discrete series of the 
Virasoro algebra.

\begin{lmm}[\cite{Mi2}]
$V_2$ decomposes into 
$$V_2=\R e\ \oplus \ E^e(0)\ \oplus \ E^e(\hf)\ \oplus \ E^e(\st),$$
where $E^e(h)$ denotes the eigenspace of $e(1)$ with eigenvalue $h$. 
\end{lmm}

\pr
As a ${\rm VA}(e)$-module, $V$ decomposes into a direct sum of copies of 
$L(\hf,0)$, $L(\hf,\hf)$ and $L(\hf,\st)$ since ${\rm VA}(e)$ is a 
rational 
VOA isomorphic to $L(\hf,0)$ 
and $L(\hf,0)$, $L(\hf,\hf)$, $L(\hf,\st)$ are the set of irreducible 
modules. 
So $V_2$ decomposes into a direct sum of 
eigenspaces of $e(1)$ with eigenvalues $0$, $2+r$, $\hf+r$ and 
$\st+r$ $(r=0,1,2,...)$. If $e(1)v=(r+h)v$ for $v\in V_2$ and 
$r+h\not=0,2,\hf,\st$, then $v$ is not in the top module of $L(c,h)$ 
as a ${\rm VA}(e)$-module and so 
there is an element $u\in {\rm VA}(e)_m$ such that 
$u(m)v\not=0$ in $L(\hf,h)$, which contradicts $u_mv\in V_1=\{0\}$. 
Since $V$ has a definite invariant bilinear form, if a central 
charge of 
a conformal vector is less than one, it should be one of minimal 
discrete 
series. Therefore $e$ is indecomposable and so $E^e(2)=\R e$, 
see \cite{Mi2}.
So we obtain the desired decomposition.
\prend

We note that from the fusion rules of $L(\hf,0)$-modules: \\
$$\begin{array}{l}
L(\hf,h)\times L(\hf,k)=L(\hf,k)\times L(\hf,h) 
\mbox{  for any }h,k=0,\hf,\st \cr
\vspace{-3mm}\cr
L(\hf,0)\times L(\hf,k)=L(\hf,k) \mbox{  for any }k=0,\hf,\st\cr
\vspace{-3mm}\cr
L(\hf,\hf)\times L(\hf,\hf)=L(\hf,0) \cr
\vspace{-3mm}\cr
L(\hf,\hf)\times L(\hf,\st)=L(\hf,\st) \hfill \mbox{   and}\cr
\vspace{-3mm}\cr
L(\hf,\st)\times L(\hf,\st)=L(\hf,0)+L(\hf,\hf) 
\end{array}$$
we have
$$\begin{array}{ll}
ab\in E^e(h) &\mbox{  for } a\in \R e+E^e(0), b\in E^e(h) 
\mbox{ and any }h \cr
\vspace{-3mm}\cr
ab\in \R e+E^e(0) &\mbox{  for } a, b\in E^e(\hf) \cr
\vspace{-3mm}\cr
ab\in E^e(\st) &\mbox{  for } a\in \R e+E^e(0)+E^e(\hf), b\in E^e(\st) 
\hfill \mbox{   and }\cr
\vspace{-3mm}\cr
ab\in \R e+E^e(0)+E^e(\hf) &\mbox{  for } a, b\in E^e(\st).  
\end{array}$$

Let $e$ and $f$ be two rational conformal vectors with central 
charge $\hf$ 
and assume that $\tau_e\tau_f$ is of order three. 
Set $G=<\tau_e,\tau_f>$. $G$ is isomorphic to a Symmetric group $S_3$ 
on three letters. 
The purpose in this section 
is to determine a subalgebra $\CG$ of $V_2$ generated by $e$ 
and $f$. \\
 
$(\tau_e\tau_f)^3=\tau_e\tau_f\tau_e\tau_f\tau_e\tau_f=1$ implies 
$\tau_e^{-1}\tau_f\tau_e=\tau_f^{-1}\tau_e\tau_f$. 
Since $\tau_{\phi(e)}=\phi \tau_e\phi^{-1}$ for any automorphism 
$\phi$ of $V$ and $\tau$ gives arise a one-to-one correspondence between 
conformal vectors and $\tau$-involutions by the assumption $(3)$ and 
the statement after $(3)$, 
we obtain  
  $$f^{\tau_e}=e^{\tau_f}.$$ 
We also note that $\tau_e^{-1}\tau_f\tau_e(e)=f$ and 
$\tau_e^{-1}\tau_f\tau_e(f)=e$. 
To simplify the 
notation, we assume that $V$ coincides with a VOA ${\rm VA}(\omega,e,f)$ 
generated by $\omega, e, f$, where $\omega$ is the 
Virasoro element of $V$.  \\

As it is well known, for a conformal vector $v$, 
its central charge is $2\langle v,v\rangle$. 
In particular, if $e$ is a conformal 
vector with central charge $\hf$, then $\langle e,e\rangle=\frac{1}{4}$. 
Set $\langle e,f\rangle =\frac{\la}{4}$. 
Using the decomposition: 
$V_2=\R e\oplus E^e(0)\oplus E^e(\hf)\oplus E^e(\st)$, we obtain 
$$  f=\la e+a+b+c    $$
where $\la\in \R$, $a\in E^e(0)$, $b\in E^e(\hf)$ and $c\in E^e(\st)$. 
By the definition of $\tau_e$, 
$$f^{\tau_e}=\la e+a+b-c. $$
Similarly, we have 
$$\begin{array}{rl}
  e=&\la f+g+h+i \qquad \mbox{  and}\cr
e^{\tau_f}=&\la f+g+h-i \cr
\end{array}$$
where $g\in E^f(0)$, $h\in E^f(\hf)$ and $i\in E^f(\st)$.

It follows from $e^{\tau_f}=f^{\tau_e}$ that $\la e+a+b-c=e-2i$ and so 
we obtain 
$$i=\frac{(1-\la)}{2}e-\hf a-\hf b+\hf c.   \eqno{(2.1)}$$ 
$ef=fe$ implies  $2\la e+\hf b+\st c=2\la f+\hf h+\st i$ 
and so we have 
$$\begin{array}{rl}
   h=&(4\la-\st)(1-\la)e+(\frac{1}{16}-4\la)a
         +(\frac{17}{16}-4\la)b+(\frac{1}{16}-4\la)c 
\end{array}\eqno{(2.2)}$$ 
On the other hand, since $e=\la f+g+h+i$ and $f=\la e+a+b+c$, we obtain 
$$\begin{array}{rl}
   g=&(1-\la)(\frac{9}{16}-3\la)e+(3\la+\frac{7}{16})a
         +(3^la-\frac{9}{16})b+(3\la-\frac{9}{16})c  
\end{array}\eqno{(2.3)}$$

By $ff=2f$, we have 
$$  2\la e+2a+2b+2c
=2\la^2 e+\la b+\frac{\la}{8}c+aa+bb+cc+2ab+2ac+2bc.   $$
Comparing the components in $E^e(\st)$, 
$\frac{\la}{8}c+2(a+b)c-2c=0$ and so we get 
$$(a+b)c=c-\frac{\la}{16}c.  \eqno{(2.4)}$$
Since $\langle e,aa\rangle=\langle ae,a\rangle=\langle 0,a\rangle=0$, 
we have 
$$ aa\in E^e(0).$$ 
Since $bb\in \R e\oplus E^e(0)$, we denote it by 
$$bb=(bb)_ee+(bb)_0, $$
where $(bb)_e\in \C, (bb)_0\in E^e(0)$. Decompose $cc$ as 
$$cc=\{(2\la-2\la^2-(bb)_e)e\}+\{2a-aa-(bb)_0\}
+\{(2-\la-2a)b\}\in \R e\oplus E^e(0)\oplus E^e(\hf). \eqno{(2.5)}$$

\noindent
From $fg=0$, 
$$0=(\la e+a+b+c)\left((\la-1)(3\la-\frac{9}{16})e+(3\la+\frac{7}{16})a+
(3\la-\frac{9}{16})b+(3\la-\frac{9}{16})c\right). $$
The components in $E^e(\st)$ are  
$$\begin{array}{l}
\dsp{0=\frac{\la}{16}(3\la-\frac{9}{16})c+(3\la-\frac{9}{16})(a+b)c
+(\la-1)(3\la-\frac{9}{16})\st c } \cr
\vspace{-3mm}\cr
\mbox{}\qquad \qquad \dsp{+(3\la+\frac{7}{16})(a+b)c-bc}
\end{array} $$
and so we have
$$ bc=\frac{23}{2^8}(2^4\la-1)c \eqno{(2.6)}$$
and
$$ ac=\frac{93}{2^8}(3-16\la)c. \eqno{(2.7)}$$

\noindent
The components in $E^e(\hf)$ are  
$$ \begin{array}{l}
\dsp{0=\frac{\la}{2}(3\la-\frac{9}{16})b+(3\la-\frac{9}{16})ab 
+\hf(\la-1)(3\la-\frac{9}{16})b }\cr
\vspace{-3mm}\cr
\mbox{}\qquad\qquad \dsp{+(3\la+\frac{7}{16})ab+
(3\la-\frac{9}{16})(cc)_{\hf}}
\end{array}$$
where $(cc)_{\hf}$ denotes the component of $cc$ in $E^e(\hf)$, 
that is, $(2-\la)b-2ab$.  
Hence 
$$\begin{array}{rl}
0=&\dsp{(\la(48\la-9)+(96\la-18)a+(48\la^2-71\la+9)+(96\la+14)a }\cr
\vspace{-3mm}\cr
&\mbox{}\qquad \dsp{+(48\la-9)(4-2\la-4a))b  }\cr
\vspace{-3mm}\cr
=&\dsp{(48\la^2-9\la+48\la^2-71\la+9+192\la-36-96\la^2+18\la }\cr
\vspace{-3mm}\cr
&\mbox{}\qquad \dsp{+192\la a-4a-192\la a+36a) }\cr
\vspace{-3mm}\cr
=&\dsp{(130\la-27+32a)b}
\end{array}$$
and so we get 
$$ ab=\frac{9}{2^5}(3-2^4\la)b.  \eqno{(2.8)}$$

\noindent
The components in $E^e(0)$ are 
$$\begin{array}{rl}
0=&\dsp{(3\la+\frac{7}{16})aa+(3\la-\frac{9}{16})(bb)_0
+(3\la-\frac{9}{16})(cc)_0 }\cr
\vspace{-3mm}\cr
=&\dsp{(3\la+\frac{7}{16})aa+(3\la-\frac{9}{16})(bb)_0
+(3\la-\frac{9}{16})(2a-aa-(bb)_0) }\cr
\vspace{-3mm}\cr
=&\dsp{3\la aa+\frac{7}{16}aa+3\la (bb)_0-\frac{9}{16} (bb)_0
+3\la 2a- 3\la aa-3\la (bb)_0}\cr
\vspace{-3mm}\cr
&\mbox{}\qquad\dsp{-2a\frac{9}{16}+\frac{9}{16}aa+\frac{9}{16}(bb)_0} 
\cr
\vspace{-3mm}\cr
=&\dsp{aa+6\la a-\frac{9}{8}a} 
\end{array}$$
and so we obtain 
$$ aa=\frac{3}{2^3}(3-2^4\la) a.   \eqno{(2.9)} $$

\noindent
Using $fh=\hf h$, we have 
$$\begin{array}{l}
\dsp{(\la e+a+b+c)
\{(\frac{65}{16}\la-4\la^2-\frac{1}{16})e+(\frac{1}{16}-4\la)a}\cr
\vspace{-3mm}\cr
\mbox{}\qquad +(\frac{17}{16}-4\la)b+(\frac{1}{16}-4\la)c\} \cr
\vspace{-3mm}\cr
\mbox{}\dsp{=\hf\{(\frac{65}{16}\la-4\la^2-\frac{1}{16})e
+(\frac{1}{16}-4\la)a
+(\frac{17}{16}-4\la)b+(\frac{1}{16}-4\la)c\}}.
\end{array}$$
Comparing the components in $W(0)$, we have
$$\begin{array}{l}
(\frac{1}{16}-4\la)aa+(\frac{17}{16}-4\la)(bb)_0
+(\frac{1}{16}-4\la)(cc)_0
=(\frac{1}{32}-2\la)a,
\end{array} $$
where $(cc)_0$ denotes the component of $cc$ at $E^e(0)$. We hence 
obtain
$$(bb)_0=\frac{3}{2^5}(2^6\la-1)a.   \eqno{(2.10)}$$

Substituting these into the expansion of $cc$, 
we have
$$cc=(2\la-2\la^2-(bb)_e)e+\frac{31}{32}a+(8\la+\frac{5}{16})b 
\eqno{(2.11)}$$

Therefore we obtain: 

\begin{lmm}
$\CG=\R e+\R a+\R b+\R c$ is a subalgebra. Furthermore, 
the symmetric group $G=<\tau_e,\tau_f>$ acts on $\CG$.
\end{lmm}

\pr
We have already shown that $\CG$ is a subalgebra. It is generated by 
$e$ and $f$ and it also contains $e^{\tau_f}$. Hence $\CG$ is invariant 
under the actions of $\tau_e$ and $\tau_f$, that is, $\CG$ is 
$G$-invariant. 
\prend

\section{Inner product}
In this section, we will calculate the inner products of elements 
in $\CG$ and show that $\la$ is either $\frac{1}{2^6}$ or 
$\frac{13}{2^8}$. 

Since  
$$
\langle ac,c\rangle =\langle a,cc\rangle 
=\frac{31}{32}\langle a,a\rangle $$
and 
$$\langle ac,c\rangle=\frac{93}{2^8}(3-2^4\la)\langle c,c\rangle,$$
we have    
$$\langle a,a\rangle=\frac{3}{8}(3-2^4\la)\langle c,c\rangle. 
\eqno{(3.1)}$$

From 
$$
\langle bc,c\rangle=\langle b,cc\rangle 
=(8\la+\frac{5}{2^4})\langle b,b\rangle \mbox{ and }
\langle bc,c\rangle=\frac{23}{2^8}(2^6\la-1)\langle c,c\rangle,$$
we obtain 
$$\langle b,b\rangle=\frac{23(2^6\la-1)}{2^4(2^7\la+5)}
\langle c,c\rangle 
\eqno{(3.2)} $$

We also have 
$$\frac{9}{32}(3-2^4\la)\langle b,b\rangle=
\langle ab,b\rangle =\langle a,bb\rangle
=\frac{3}{2^5}(2^6\la-1)\langle a,a\rangle. \eqno{(3.3)} $$
Substituting (3.1) and (3.2) into (3.3),  
$$(2^6\la-1)(-2^4\la)\langle c,c\rangle
=\frac{23(2^6\la-1)}{2(2^7\la+5)}\langle c,c\rangle. \eqno{(3.4)}$$
If $\langle c,c\rangle=0$, then $c=0$ and so 
$f^{\tau_e}=f$ and $e=f^{\tau_e\tau_f}=f$, 
which contradicts $e\not=f$. Hence $\langle c,c\rangle\not=0$ and 
$\la$ is one of 
$$\frac{3}{2^4}, \frac{1}{2^6}, \mbox{ or } \frac{13}{2^8}. 
\eqno{(3.5)}$$

On the other hand, we obtain
$$\langle c,c\rangle=16\langle ec,c\rangle
=16\langle e,cc\rangle
=4(2\la-2\la^2-(bb)_e). \eqno{(3.6)}$$

Since 
$$(\tau_e^{-1}\tau_f^{-1}\tau_e)\tau_e(\tau_e^{-1}\tau_f\tau_e)
=\tau_e\tau_f\tau_e\tau_f\tau_e=\tau_f, $$
$\tau_e^{-1}\tau_f\tau_e(e)=f$ and $\tau_e^{-1}\tau_f\tau_e(f)=e$. 
Since $c$ and $i$ are 
uniquely defined by $e$ and $f$, $\tau_e\tau_f\tau_e(c)=i$ and 
so we obtain 
$\langle c,c\rangle=\langle i,i\rangle$. Therefore, 
$$\begin{array}{rl}
\dsp{\langle c,c\rangle=}&\dsp{\langle i,i\rangle
=\frac{1}{4}\langle (1-\la)e-a-b+c, (1-\la)e-a-b+c\rangle }\cr 
\vspace{-3mm}&\cr
=&\dsp{\frac{1}{4}\{(1-\la)^2\frac{1}{4}
+\langle a,a\rangle +\langle b,b\rangle 
+\langle c,c\rangle\} }\cr
\vspace{-3mm}&\cr
=&\dsp{\frac{1}{16}(1-2\la)+\frac{1}{16}=\frac{1}{8}(1-\la)}
\end{array}\eqno{(3.7)}$$
and so we have 
$$\langle a,a\rangle =\frac{3(1-\la)(3-2^4\la)}{2^6} \eqno{(3.8)} $$
$$\langle b,b\rangle =\frac{23}{2^7}\frac{(1-\la)(2^6\la-1)}{2^7\la+5} 
\qquad  \mbox{ and } \eqno{(3.9)} $$
$$\langle c,c\rangle =\frac{1-\la}{8} \eqno{(3.10)} $$

Using $\langle f,f\rangle =\frac{1}{4}$, we have 
$$
\begin{array}{rl}
\frac{1}{4}
=&\dsp{\la^2\frac{1}{4}+\frac{-3\cdot 2^4\la+9}{8}\langle c,c\rangle 
+\frac{23(2^6\la-1)}{(2^7\la+5)2^4}\langle c,c\rangle 
+\langle c,c\rangle } \cr
=&\dsp{\la^2\frac{1}{4}+\frac{-3\cdot 2^4\la+9}{8}\hf(1-\la)
+\frac{23(2^6\la-1)}{(2^7\la+5)2^4}\hf(1-\la)+\hf(1-\la)} 
\end{array}$$
and 
$$\begin{array}{rl}
0=&\dsp{(\la^2-1)\frac{1}{4}+\frac{-3\cdot 2^4\la+9}{8}\hf(1-\la)
+\frac{23(2^6\la-1)}{(2^7\la+5)2^4}\hf(1-\la)+\hf(1-\la) }\cr
=&\dsp{(1-\la)\{
(-\la^2-1)\frac{1}{4}+\frac{-3\cdot 2^4\la+9}{8}\hf
+\frac{23(2^6\la-1)}{(2^7\la+5)2^4}\hf+\hf\} }\cr
=&\dsp{\frac{(1-\la)(2^6\la-1)(13-2^8\la)}{2^4(2^7\la+5)}}.
\end{array}$$
It follows from (3.5) that $\la$ is either 
$\dsp{\frac{1}{2^6}}$ or $\dsp{\frac{13}{2^8}}$. 

As we mentioned in the introduction, 
the both cases occur in the monstrous Griess algebra. 
So we have the following theorem. 

\begin{thm} 
Let $V$ be a VOA of moonshine type satisfying the assumption $(3)$ and 
let 
$e$ and $f$ be two distinct rational conformal vectors with 
central charge $\hf$.  If $\tau_e\tau_f$ is of order three, then 
$\langle e,f\rangle $ is either $\frac{1}{2^8}$ or $\frac{13}{2^{10}}$. 
\end{thm}

We next study a structure of a subVA ${\rm VA}(e,f)$ generated by 
$e$ and $f$. First we will show that 
${\rm VA}(e,f)$ has a Virasoro element and $\CG$ contains it. 

Since $aa=\frac{3}{8}(3-16\la)a$, 
$$\omega_1=\frac{16}{(9-48\la)}a \eqno{(3.11)}$$ 
is a conformal vector. 

\begin{lmm}
$e+\omega_1$ is a Virasoro element of ${\rm VA}(e,f)$. 
\end{lmm}

\pr
By direct calculation using (2.7),(2.8),(2.9), 
$$(e+\omega_1)e=2e, \quad (e+\omega_1)a=2a, \quad (e+\omega_1)b
=2b, \quad (e+\omega_1)c=2c. $$
So it is sufficieint to show 
that $(e+\omega_1)(0)$ satisfies the derivation property 
for $v\in {\rm VA}(e,f)$, that is, $Y((e+\omega_1)(0)v,z)
=\frac{d}{dz}Y(v,z)$. 
As a ${\rm VA}(e)\cong L(\hf,0)$-module, $a$ is a highest 
weight vector in 
$L(\hf,0)$ and so $e(0)a=0$. By the skew-symmetry property, 
$a(0)e=0$ and so 
$$Y((e+\omega_1)(0)e,z)=Y(e(0)e,z)=\frac{d}{dz}Y(e,z).$$ 
Since $e+\omega_1$ is a unique element of $\CG$ satisfying 
$(e+\omega_1)v=2v$ for all $v\in \CG$, $G$ fixes it and so 
$Y((e+\omega_1)(0)f,z)=\frac{d}{dz}Y(f,z)$. It is shown in 
\cite{Li} that 
if $a(z)$ and $b(z)$ satisfy $D$-derivation property, 
then so does $a(z)_nb(z)$ for any integer $n$, where $a(z)_nb(z)$ 
denotes $n$-th product. 
Hence $(e+\omega_1)(0)$ has the derivation property for all 
elements in ${\rm VA}(e,f)$. 
\prend

From now on, we assume 
that $V={\rm VA}(e,f)$ because it has a Virasoro element. 
Since it is convenient to use a linear representation of 
$<\tau_e\tau_f>$, we will treat 
the complexification $\C\otimes_{\R} V$ of $V$.  So from now on, 
$V$ denotes $\C\otimes_{\R} V$.

\section{The case $\la=\frac{1}{2^6}$}
First we will study the case $\la=\frac{1}{2^6}$. 
In this case, $\langle b,b\rangle=0$ and so $b=0$. Namely, 
$\CG$ is of dimension three. 
We will show $({\rm VA}(e,f))_2=\CG$. 

The structure of $\CG$ is given by 
$$\begin{array}{l}
\dsp{aa=\frac{33}{32}a }\cr
\vspace{-3mm}\cr
\dsp{ac=\frac{7\times 11\times 13}{2^{10}}c  }\cr
\vspace{-3mm}\cr
\dsp{cc=\frac{3^2\times 7}{2^{11}}e+\frac{31}{32}a+\frac{7}{2^4}b. }
\end{array}$$

In $\CG=\C e+\C a+\C c$, a Virasoro element decomposes into 
an orthogonal sum of two conformal vectors $e$ with central charge 
$\hf$ and $\omega_1$ with central charge $\frac{21}{22}$. 
The both 
${\rm VA}(e)\cong L(\hf,0)$ and ${\rm VA}(\omega_1)\cong 
L(\frac{21}{22},0)$ 
are minimal series and rational. In particular, 
$\omega_1$ is indecomposable and so $E^e(0)=\C a$ by \cite{Mi2}.  \\

\begin{lmm}
$E^e(\hf)=0$.
\end{lmm}

\pr
Viewing $V$ as a $L(\hf,0)\otimes L(\frac{21}{22},0)$-module, 
$V$ is a direct sum of irreducible 
$L(\hf,0)\otimes L(\frac{21}{22},0)$-modules. 
If $E^e(\hf)\not=0$, then $L(\frac{21}{22},0)$ has to have an 
irreducible 
module with a highest weight $\frac{3}{2}$. 
However, there is no $\frac{3}{2}$ in the list 
of highest weights of $L(\frac{21}{22},0)$-modules: 
$$  \frac{(12r-11s)^2-1}{4\times 11\times 12}  
\qquad 1\leq s\leq r \leq 10. $$
\prend

Actually, the highest weights of irreducible 
$L(\frac{21}{22},0)$-modules are 
$$\begin{array}{l}
\dsp{0, \  \frac{3}{16\times 11}, \  \frac{1}{16}, \  
\frac{5}{2^4\times 11}, \  \frac{1}{22}, \  
\frac{35}{2^4\times 3\times 11}, \  \frac{1}{11}, \  
\frac{21}{2^4\times 11}, \  
\frac{5}{33}, \  \frac{3}{16}, \  
\frac{7}{22}, \  }\cr
\vspace{-3mm}\cr
\dsp{\frac{65}{2^4\times 11}, \  \frac{14}{33}, \  
\frac{85}{2^4\times 11}, \  
\frac{6}{11}, \  
\frac{323}{2^4\times 3\times 11}, \  \frac{15}{22}, \  
\frac{133}{2^4\times 11}, \  
\frac{5}{6}, \  \frac{13}{11}, \ 
\frac{225}{2^4\times 11}, \  \frac{91}{2\times 3\times 11}, \  }
\cr \vspace{-3mm}\cr
\dsp{\frac{261}{2^4\times 11}, \  \frac{35}{22}, \  
\frac{899}{2^4\times 3\times 11}, \  
\frac{20}{11}, \  \frac{31}{2^4}, \  \frac{57}{11}, \  
\frac{481}{2^4\times 11}, \  
\frac{95}{33}, \  
\frac{533}{2^4\times 11}, \  \frac{35}{11}, \  
\frac{1763}{2^4\times 3\times 11}, \  }\cr
\vspace{-3mm}\cr
\dsp{\frac{7}{2}, \  \frac{50}{11}, \  
\frac{833}{2^4\times 11}, \  \frac{325}{66}, \  
\frac{901}{2^4\times 11}, \  
\frac{117}{22}, \  \frac{265}{2^4\times 3}, \  
\frac{155}{11}, \  \frac{1281}{2^4\times 11}, \  
\frac{248}{33}, \  
\frac{1365}{2^4\times 11}, \  }\cr
\vspace{-3mm}\cr
\dsp{8, \  \frac{121}{12}, \  \frac{1825}{2^4\times 11}, \  
\frac{703}{66}, \  
\frac{175}{2^4}, \  \frac{301}{22}, \  
\frac{2465}{2^4\times 11}, \  \frac{43}{3}, \  \frac{196}{11}, \  
\frac{291}{16}, \  \frac{45}{2}}
\end{array}$$ 

Since $V$ has integer weights, $V$ is a direct sum of copies of 
$$ \begin{array}{l}
L(\hf,0)\otimes L(\frac{21}{22},0), \quad
L(\hf,0)\otimes L(\frac{21}{22},8), \quad 
L(\hf,\hf)\otimes L(\frac{21}{22},\frac{7}{2}), \cr
\vspace{-3mm}\cr
L(\hf,\hf)\otimes L(\frac{21}{22},\frac{45}{2}), \quad
L(\hf,\st)\otimes L(\frac{21}{22},\frac{31}{16}), \quad 
L(\hf,\st)\otimes L(\frac{21}{22},\frac{175}{16})
\end{array}$$
as a $L(\hf,0)\otimes L(\frac{21}{22},0)$-module.

\begin{lmm}
$V_2=\CG$ and $\dim V_2=3$. 
\end{lmm}

\pr
Suppose false and set $T=\CG^{\perp}$. 
Since $\dim E^e(0)=1$ and $E^e(\hf)=0$, 
$T\subseteq E^e(\st)$ and so 
$ev=\frac{1}{16}v$ and $\tau_e(v)=-v$ for $v\in T$.
Since $G$ acts on $T$. 
$fv=(\frac{1}{2^6}e+a+c)v=\st v$, which implies $cv=0$ and 
$av=\frac{2^6-1}{2^{10}}$. 
Moreover, $\tau_e\tau_f$ acts $T$ as $1$. In particular, 
$vu\in (\C e+E^e(0))^{G}=\C \omega$ for $u,v\in T$. 
However, since 
$0=\langle a,uu\rangle=\langle au,u\rangle
=\frac{2^6-1}{2^{10}}\langle u,u\rangle$, we have a 
contradiction. 
\prend

\begin{thm} 
Assume that $\tau_e\tau_f$ is of order three and 
$\langle e,f\rangle=\frac{1}{2^8}$. Then 
${\rm VA}(e,f)$ is a VOA with central charge $\frac{16}{11}$ and 
a Griess algebra $({\rm VA}(e,f))_2$ is of dimension 3. 
${\rm VA}(e,f)$ contains $L(\hf,0)\otimes L(\frac{21}{22},0)$. 
\end{thm}

\section{The case $\la=\frac{13}{2^8}$}
We are very interested in this case because each  
edge in the diagram of $Y_{5,5,3}$ or the graph of $21$ involutions is 
corresponding to this case. 

As we showed in the previous section, the structure of $\CG$ is given by 
$$\begin{array}{l}
\dsp{ab=\frac{3^2\times 5\times 7}{2^9}b }\cr
\vspace{-3mm}\cr
\dsp{aa=\frac{3\times 5\times 7}{2^7}a }\cr
\vspace{-3mm}\cr
\dsp{bb=\frac{3^9}{2^{15}}e+\frac{3^3}{2^7})a } \cr
\vspace{-3mm}\cr
\dsp{ac=\frac{5\times 7^2\times 13}{2^{12}}c } \cr
\vspace{-3mm}\cr
\dsp{bc=\frac{3^2\times 23}{2^{10}}c }\cr
\vspace{-3mm}\cr
\dsp{cc=\frac{3^5}{2^{13}}e+\frac{31}{32}a+\frac{23}{2^5}b.}
\end{array}$$

In particular, 
$$\omega_1=\frac{2^8}{105}a$$ 
is a conformal vector with central charge $\dsp{\frac{81}{70}}$.

\begin{lmm}
$\omega_1$ is indecomposable. In particular, $E^e(0)=\C \omega_1$. 
\end{lmm}

\pr
Suppose that $\omega_1$ is a sum $\omega'+\omega''$ of orthogonal 
conformal 
vectors. Since $c.c(\omega'), c.c(\omega'')\geq \hf$, 
$c.c(\omega'), c.c(\omega'')\leq \frac{23}{35}$, 
where $c.c(\omega')$ denotes a central charge of $\omega'$. 
In particular, 
${\rm VA}(\omega')$ and ${\rm VA}(\omega'')$ belong to minimal 
discrete series. However, the central charge of minimal discrete series 
which is less than 
or equal to $\frac{23}{35}$ are only $\{\hf, \frac{7}{10}\}$. 
However, $\frac{81}{70}$ is not a sum of them, 
which is a contradiction. 
\prend

\subsection{$\theta$-fixed point space}
In this subsection, we will study a $\theta$-fixed point space and 
find conformal vectors of central charge $\ff$ and $\ssv$. 
Using these conformal vectors, we will also show that 
$({\rm VA}(e,f))_2=\CG$. 

Set $\theta=\tau_e\tau_f$, which is an automorphism of $V$ of order 
three. 
We note that $e^{\theta}=f$ and $f^{\theta}=f^{\tau_e}$.
Hence $\alpha=e+f+f^{\tau_e}=(\frac{13}{2^7}+1)e+2a+2b$ 
and $\omega$ are in $\CG^{<\theta>}$, 
where $\CG^{<\theta>}=\{v\in \CG |v^{\theta}=v \}$.

By the direct calculation, we have
$$\begin{array}{l}
\dsp{ef=\frac{13}{2^7}e+\hf b+\st c }\cr
\dsp{ff^{\tau_e}=(\frac{13}{2^7}
+\frac{13}{2^{12}}-\st)e+\st a+\frac{9}{16}b }\cr
\dsp{f^{\tau_e}e=\frac{13}{2^7}e+\hf b-\st c} 
\end{array} \eqno{(5.1)}$$ 
Hence 
$$\alpha\alpha
=2e+2f+2f^{\tau_e}+(\frac{39}{2^6}+\frac{13}{2^{11}}-\frac{1}{8})e
+\frac{1}{8}a+\frac{25}{16}b $$
and so 
$$ \alpha\alpha
=\frac{57}{16}\al+(\frac{9\times 13}{2^8}-\frac{27}{16})\omega. 
\eqno{(5.2)}$$

Setting $\be=\frac{2^4}{3}\alpha$, we obtain 
$$\begin{array}{l}
\dsp{\be\be=19\be-35\omega }\cr
\vspace{-3mm}\cr
\dsp{\langle \be,\be\rangle=\frac{47}{2} }\cr
\vspace{-3mm}\cr
\dsp{\langle \be,\omega\rangle=4 \qquad \mbox{  and }}\cr
\vspace{-3mm}\cr
\dsp{\langle \omega,\omega\rangle=\frac{29}{35}.}
\end{array}$$

It follows from a direct calculation that 
$$\omega_2=\frac{2(7\omega-\be)}{9}$$ 
is a conformal vector with central charge $\frac{4}{5}$ 
and $\omega_3=\omega-\omega_2$ is a conformal vector with 
central charge $\frac{6}{7}$. We note that 
the both belong to minimal discrete series. 
In particular, we have:

\begin{prn}
$V$ contains a rational VOA 
$L(\frac{4}{5},0)\otimes L(\frac{6}{7},0)$ whose Virasoro 
element is equal to $\omega$. 
\end{prn}

We next view $V$ as a $L(\ff,0)\otimes L(\ssv,0)$-module.
Let $\ep$ be a cubic root of unity and set 
$\ga=e+\ep f+\ep^2 f^{\tau_e}$. Namely, $\ga^{\theta}=\ep^{-1}\ga$. 
Then we have 
$$\begin{array}{rl}
\alpha\ga=&(e+f+f^{\tau_e})(e+\ep f+\ep^2 f^{\tau_e}) 
\cr
\vspace{-3mm}&\cr
=&2(e+\ep f+\ep^2 f^{\tau_e})+(\ep+1)ef+(\ep^2+1)ef^{\tau_e}
+(\ep+\ep^2)ff^{\tau_e} \cr
\vspace{-3mm}&\cr
=&2\be+\st \ga=\frac{33}{16}\ga,  
\end{array}$$
which implies 
$$\be\ga=(\frac{2^4}{3}\alpha)\ga=11\ga $$
and 
$$\omega_2\ga=\tth \ga. $$
We note $\{e, f, f^{\tau_e}\}$ is a linearly independent set and 
so $\ga\not=0$. $\theta=\tau_e\tau_f$ acts on $\C\ga$ as $\ep^{-1}$.   
Similarly, $\omega_2\overline{\ga}=\tth \overline{\ga}$. 
Therefore, 
$V$ contains $U^1\oplus U^2$ where 
$U^i\cong L(\ff,\tth )\otimes L(\frac{6}{7},\frac{4}{3})$ and 
$\theta=\tau_e\tau_f$ acts on $U^i$ as $\ep^i$ and $\tau_e$ exchange 
$U^1$ and $U^2$. \\

Therefore,  
$$e, f\in L(\ff,0)\otimes L(\frac{6}{7},0)
+L(\ff,\tth )\otimes L(\frac{6}{7},\frac{4}{3})
+L(\ff,\tth )\otimes L(\frac{6}{7},\frac{4}{3}).$$
From fusion rules of $L(\ff,0)$-modules:
$$\begin{array}{l}
\dsp{L(\ff,\tth )\times L(\ff,\tth )=L(\ff,0)+L(\ff,\tth )+L(\ff,3) }
\cr
\vspace{-3mm}\cr
\dsp{L(\ff,3)\times L(\ff,\tth )=L(\ff,\tth ) }\cr
\vspace{-3mm}\cr
\dsp{L(\ff,3)\times L(\ff,3)=L(\ff,0)},
\end{array}$$
we obtain that $V={\rm VA}(e,f)$ is a direct sum of copies of 
$L(\ff,0)$, $L(\ff,\tth )$ and $L(\ff,3)$ 
as a $L(\ff,0)$-module, see Appendix.

The set of highest weights of irreducible $L(\frac{6}{7},0)$-modules are 
$$  \frac{(7r-6s)^2-1}{4\times 6\times 7} \quad 1\leq s\leq r\leq 5, $$
that is, 
$$ \frac{1}{56}, \frac{1}{21}, \frac{5}{56},\frac{1}{7}, 
\frac{3}{8}, \frac{10}{21}, \frac{33}{56}, \frac{5}{7}, \frac{4}{3}, 
\frac{85}{56}, \frac{12}{7}, \frac{23}{8}, \frac{22}{7}, \mbox{  and }5.
$$
Since the weights in $V$ are integers, 
the possibilities of $L(\ff,0)\otimes L(\ssv,0)$-submodules are 
$$ \begin{array}{l}
\dsp{L(\ff,0)\otimes L(\ssv,0), \quad L(\ff,0)\otimes L(\ssv,5), }
\cr
\vspace{-3mm}\cr
\dsp{L(\ff,\tth )\otimes L(\ssv,\frac{4}{3}), }\cr
\vspace{-3mm}\cr
\dsp{L(\ff,3)\otimes L(\ssv,0), \quad L(\ff,3)\otimes L(\ssv,5)}.
\end{array}$$

The fusion rules for $L(\ssv,0)$-modules $L(\ssv,\frac{4}{3})$ 
and $L(\ssv,5)$ are 
$$\begin{array}{l}
\dsp{L(\ssv,5)\times L(\ssv,5)=L(\ssv,0) }\cr
\vspace{-3mm}\cr
\dsp{L(\ssv,5)\times L(\ssv,\frac{4}{3})=L(\ssv,\frac{4}{3}) }\cr
\vspace{-3mm}\cr
\dsp{L(\ssv,\frac{4}{3})\times L(\ssv,\frac{4}{3})=
L(\ssv,0)+L(\ssv,\frac{4}{3})+L(\ssv,5)}, 
\end{array}$$
see \cite{W}.

Our last purposes is to show the following theorem. 

\begin{thm}
$({\rm VA}(e,f))_2=\CG$. 
\end{thm}

\pr
Assume $\CG\not=V_2$ and set $T=\CG^{\perp}$. 
First we will show the following lemma. 

\begin{lmm}
$T^{<\theta>}=0$ and $(V_2)^{<\theta>}=\C \omega_2+\C \omega_3$. 
\end{lmm}

\pr
Clearly, $V$ contains only one copy of 
$L(\ff,0)\otimes L(\ssv,0)$. So the eigen value of $\omega_2$ on $T$ 
is only $\tth $.
It is clear that $G$ acts on $T$. 
Suppose $P(\hf)=T^{\theta}\cap E^e(\hf)\not=0$. Namely, 
$\tau_e, \tau_f, \tau_{f^{\tau_e}}$ act $P(\hf)$ as $1$. 
In this case, for $u\in P(\hf)$, 
$eu=\hf u$, $fu=\hf u$, $f^{\tau_e}u=\hf u$. 
Hence $cu=0$. Moreover, $bu\in (\C e+E^e(0))\cap T=\{0\}$. 
Since 
$$2u=\omega u=(e+\omega_1)u=(e+\frac{2^8a}{105})u
=\hf u+\frac{2^8}{105}au,$$ 
we obtain 
$\dsp{au=\frac{315}{2^9}u}$. However, we have
$$\hf u=fu=(\frac{13}{2^8}e+a+b+c)u=\frac{13}{2^9}u+\frac{315}{2^9}u
=\frac{328}{512}u,$$ 
a contradiction. 
We hence have $P(\hf)=0$. 
We next assume $P(\st)=T^{\theta}\cap E^e(\st)\not=0$. Namely, 
$\tau_e, \tau_f, \tau_{f^{\tau_e}}$ act on 
on $P(\st)$ as $-1$. In this case, for $u\in P(\st)$, 
$eu=\st u$, $fu=\st u$, $f^{\tau_e}u=\st u$ by the 
definition of $\tau$-automorphism. 
Hence $cu=\hf(f-f^{\tau_e})u=0$ and $\al u=\{(2\la+1)e+2a+2b\}u
=eu+2(\la e+a+b)u=\frac{3}{16}u$. Since 
$\dsp{\omega_2=\frac{2(21\omega-16\alpha)}{27}}$, we have 
$\dsp{\tth u=\omega_2u=\frac{2(\be-7\omega)}{-9}u
=\frac{2(1-14)}{-9}u}$, 
a contradiction. Therefore $T^{\theta}=0$ and it is clear that 
$\CG^{<\theta>}=\C \omega_2+\C \omega_3$. 
\prend

We next view $V$ as a $G$-module. 
The symmetric group $G$ on three letters has 
$3$ irreducible modules $\C(+)$, $\C(-)$ and $D$, where $\C(+)$ is a 
trivial module, $\theta$ and $\tau_e$ act on $\C(-)$ as $1$ and $-1$, 
respectively, 
and $D$ is an irreducible module of dimension two. 
Therefore, $V$ has a subVOA $V^{G}$ such that 
$V$ decomposes into the direct sum 
$$ V=V^{G}\oplus U\otimes \C(-) \oplus W\otimes D,  $$
where $U$ and $W$ are irreducible $V^{G}$-modules, 
see \cite{DM}.

Since $\omega_2, \omega_3\in {\rm VA}(e,f)^{G}$, 
$V^{G}$ contains $L(\ff,0)\otimes L(\ssv,0)$. 

Recall
$\ga=e+\ep f+\ep^2 f^{\tau_e}$ and $\theta=\tau_e\tau_f$ acts $\ga$ 
as $\ep^{-1}$. It is easy to check that $\ga\not=0$. 
Therefore $W$ contains 
$L(\ff,\frac{2}{3})\otimes L(\ssv,\frac{4}{3})$. 
Lemma 5.3 implies the highest weight of $U$ is greater than 3. \\

\noindent
Case 1. \\
First assume $V^{G}=L(\ff,0)\otimes L(\ssv,0)$. 
Then $W=L(\ff,\frac{2}{3})\otimes L(\ssv,\frac{4}{3})$ and 
so $V_2=(L(\ff,0)\otimes L(\ssv,0))_2\oplus (W\otimes D)_2$ is of  
dimension four. Namely we have the desired result. \\

\noindent
Case 2. \\
Assume $V^{G}\not=L(\ff,0)\otimes L(\ssv,0)$. 
Then 
$$V^{<\theta>}=(\ff,0)\otimes L(\ssv,0)\oplus 
L(\ff,0)\otimes L(\ssv,5)\oplus L(\ff,3)\otimes L(\ssv,0)\oplus 
L(\ff,3)\otimes L(\ssv,5),$$
which is a subVOA.  
In particular, $L(\ssv,0)\oplus L(\ssv,5)$ has a VOA structure. 

$L(\ff,0)\oplus L(\ff,3)$ is a rational VOA 
and its irreducible modules are classified in \cite{KMiY}. 
Namely, it has exactly the following six irreducible modules: 
$$ W(0), W({2\over 5}), 
W({2\over 3},+), W({1\over 15},+), W({2\over 3},-), 
W({1\over 15},-).$$
Here $h$ in $W(h)$ and $W(h,\pm)$ denotes the lowest degree 
and $W(k,-)$ is the contragredient (dual) module of $W(k,+)$ 
for $k={2\over 3}, {1\over 15}$.
As $L(\ff,0)$-modules, 
$$ \begin{array}{l}
W(0)\cong L(\ff,0)\oplus L(\ff,3), \cr
\vspace{-3mm}\cr
W({2\over 5})\cong L(\ff,{2\over 5})\oplus L(\ff,{7\over 5}),  
\cr
\vspace{-3mm}\cr
W({2\over 3},+)\cong L(\ff,{2\over 3}),  \cr
\vspace{-3mm}\cr
W({2\over 3},-)\cong L(\ff,{2\over 3}),  \cr
\vspace{-3mm}\cr
W({1\over 15},+)\cong L(\ff,{1\over 15}),  \cr
\vspace{-3mm}\cr
W({1\over 15},-)\cong L(\ff,{1\over 15}). 
\end{array} $$

For an irreducible $L(\ssv,0)\oplus L(\ssv,5)$-module, 
we will need the following lemma. 

\begin{lmm}  If $L(\ssv,0)\oplus L(\ssv,5)$ is a VOA and $X$ is an 
irreducible module contains $L(\ssv, \frac{4}{3})$, then 
$X\cong L(\ssv, \frac{4}{3})$ as $L(\ssv,0)$-modules.  
\end{lmm}

\pr 
Suppose false and let $Y^X$ be a vertex operator of 
$L(\ssv,0)\oplus L(\ssv,5)$ on $X$. 
Since $L(\ssv,0)$ is rational, $X$ is a direct sum of 
$L(\ssv,0)$-modules. 
By the fusion rule 
$L(\ssv,5)\times L(\ssv,\frac{4}{3})=L(\ssv,\frac{4}{3})$, we have 
$$X=L(\ssv,\frac{4}{3})\oplus L(\ssv,\frac{4}{3}).$$ 

Choose $0\not=I\in I{L(\ssv,\frac{4}{3})\choose L(\ssv,5)
\quad L(\ssv,\frac{4}{3})}$.  
Since $\dim I{L(\ssv,\frac{4}{3})\choose L(\ssv,5)
\quad L(\ssv,\frac{4}{3})}=1$, 
there is $\mu\in \C$ such that 
$Y^X(u,z)=\pmatrix{0 & I(u,z) \cr \mu I(u,z) &0}$ for $u\in L(\ssv,5)$ 
by choosing a suitable basis. Replacing $I$ by $\sqrt{mu}I$, we have 
$Y^X(u,z)=\pmatrix{0 & I(u,z) \cr I(u,z) &0}$.
Then $\{(w,w)\in L(\ssv,\frac{4}{3})\oplus L(\ssv,\frac{4}{3})\}$ is a 
$L(\ssv,0)\oplus L(\ssv,5)$-submodule, which contradicts 
the assumption on $X$. 
\prend

Now we go back to the proof of Theorem 5.3.  
Viewing $V$ as a $<\theta>$-module, \\
$V=V^{0}\oplus V^{1}\oplus V^2$, where $\theta$ acts 
on $V^i$ as $(\ep)^i$. Since 
$$V^0\cong (L(\ff,0)\oplus L(\ff,3))\otimes 
(L(\ssv,0)\oplus L(\ssv,5)),$$ 
we have 
$$V^i\cong L(\ff,\frac{2}{3})\otimes L(\ssv,\frac{4}{3}) $$
and so 
$\dim V_2=4$ and $V_2=\CG$. 
This completes the proof of Theorem 5.3. 
\prend

On $X=L(\ssv,\frac{4}{3})$, there are two 
$L(\ssv,0)\oplus L(\ssv,5)$-module structures.  
Namely, if $Y^X(v,z)$ is a vertex operator of 
$v\in L(\ssv,0)\oplus L(\ssv,5)$, then the other is defined by 
$Y^X(v,z)$ for $v\in L(\ssv,0)$ and $-Y^X(v,z)$ for $v\in L(\ssv,5)$. 
We denote them by $L(\ssv,\frac{4}{3})^{\pm 1}$.

\begin{thm}
If $|\tau_e\tau_f|=3$ and $\langle e,f\rangle=\frac{13}{2^{10}}$, then  
${\rm VA}(e,f)$ is isomorphic to one of following: 
$$\begin{array}{ll}
(1)& (0,0)\oplus(3,0)\oplus 
\left(W(\frac{2}{3},+)\otimes L(\ssv,\frac{4}{3})\right)
\oplus \left(W(\frac{2}{3},-)\otimes L(\ssv,\frac{4}{3})\right) 
\cr 
\vspace{-3mm}\cr
(2)& (0,0)\oplus(0,5)\oplus \left(
L(\ff,\frac{2}{3})\otimes L(\ssv,\frac{4}{3})^{+1}\right)
\oplus \left(L(\ff,\frac{2}{3})\otimes L(\ssv,\frac{4}{3})^{-1}\right) 
\cr 
\vspace{-3mm}\cr
(3)& (0,0)\oplus(3,5)\oplus \left(W(\frac{2}{3},+)
\otimes L(\ssv,\frac{4}{3})^{+1}\right)
\oplus \left(W(\frac{2}{3},-)\otimes L(\ssv,\frac{4}{3})^{+1}\right) 
\cr 
\vspace{-3mm}\cr
(4)& (0,0)\oplus(3,0)\oplus(0,5)\oplus(3,5)
\oplus \left(W(\frac{2}{3},\pm)\otimes L(\ssv,\frac{4}{3})^{+1}\right)
\oplus \left(W(\frac{2}{3},\mp)\otimes L(\ssv,\frac{4}{3})^{-1}\right) 
\cr
&\mbox{in (4)}, \tau_e\tau_f \mbox{ is equal to an automorphism 
determined by 
a 3-State Potts model }\cr 
&L(\ff,0)\oplus L(\ff,3) \mbox{ defined in \cite{Mi3}}
\end{array}$$
where $(h,k)$ denotes $L(\ff,h)\otimes L(\ssv,k)$. 
We will call them VOAs of involution type $A_2$.  
\end{thm}

At last we expect the following conjecture. 

\begin{cnj}
Let $e$ and $f$ be two distinct rational conformal vectors with 
central charge $\hf$.  If $\dsp{\langle e,f\rangle =\frac{13}{2^8}}$, 
then 
$\tau_e\tau_f$ is of order three.
\end{cnj}

\section{Appendix}
\begin{center}
{Fusion rule of $L(\frac{4}{5},0)$-modules}
\end{center}
\small
$$ \begin{array}{|c||c|c|c|c|c|c|c|c|c|}
\hline
0& \frac{2}{5}& \frac{1}{40}& \frac{7}{5}& \frac{21}{40}& 
\frac{1}{15}& 3& 
\frac{13}{8}& \tth & \frac{1}{8}  \cr
\hline 
\hline
\frac{2}{5}& 0\!:\!\frac{7}{5}& \frac{1}{8}\!:\!\frac{21}{40}
& \frac{2}{5}\!:\!3& \frac{1}{40}\!:\!\frac{13}{8}
& \frac{1}{15}\!:\!\tth & \frac{7}{5}& \frac{21}{40}
& \frac{1}{15}& \frac{1}{40}  \cr
\hline 
\frac{1}{40}& \frac{1}{8}\!:\!\frac{21}{40}
& 0\!:\!\frac{7}{5}\!:\!\tth \!:\!\frac{1}{15}
& \frac{1}{40}\!:\!\frac{13}{8}
& \frac{2}{5}\!:\!3\!:\!\frac{1}{15}\!:\!\tth & 
\frac{1}{40}\!:\!\frac{13}{8}\!:\!\frac{21}{40}\!:\!\frac{1}{8}
& \frac{21}{40}& \frac{7}{5}\!:\!\frac{1}{15}
& \frac{21}{40}\!:\!\frac{1}{40}& \frac{1}{15}\!:\!\frac{2}{5}  
\cr
\hline 
\frac{7}{5}& \frac{2}{5}\!:\!3& \frac{1}{40}\!:\!\frac{13}{8}
& 0\!:\!\frac{7}{5}& \frac{1}{8}\!:\!\frac{21}{40}
& \tth \!:\!\frac{1}{15}& \frac{2}{5}
& \frac{1}{40}& \frac{1}{15}& \frac{21}{40} \cr
\hline 
\frac{21}{40}& \frac{1}{40}\!:\!\frac{13}{8}
& \frac{2}{5}\!:\!3\!:\!\frac{1}{15}\!:\!\tth 
& \frac{1}{8}\!:\!\frac{21}{40}
& 0\!:\!\frac{7}{5}\!:\!\tth \!:\!\frac{1}{15}
& \frac{1}{8}\!:\!\frac{21}{40}
\!:\!\frac{13}{8}\!:\!\frac{1}{40}& \frac{1}{40}
& \frac{2}{5}\!:\!\frac{1}{15}& \frac{1}{40}\!:\!\frac{21}{40}
& \frac{1}{15}\!:\!\frac{7}{5}  \cr
\hline 
\frac{1}{15}& \frac{1}{15}\!:\!\tth 
& \frac{1}{40}\!:\!\frac{13}{8}\!:\!\frac{21}{40}\!:\!\frac{1}{8}
& \tth \!:\!\frac{1}{15}
& \frac{1}{8}\!:\!\frac{21}{40}\!:\!\frac{13}{8}\!:\!\frac{1}{40}& 
0\!:\!\frac{7}{5}\!:\!\tth \!:\!\frac{1}{15}\!:\!3\!:\!\frac{2}{5}
& \frac{1}{15}& \frac{1}{40}\!:\!\frac{21}{40}
& \frac{2}{5}\!:\!\frac{1}{15}\!:\!\frac{7}{5}
& \frac{1}{40}\!:\!\frac{21}{40}  \cr
\hline 
3& \frac{7}{5}& \frac{21}{40}& \frac{2}{5}& \frac{1}{40}& \frac{1}{15}
& 0& \frac{1}{8}& \tth & \frac{13}{8}  \cr
\hline 
\frac{13}{8}& \frac{21}{40}& \frac{7}{5}\!:\!\frac{1}{15}& \frac{1}{40}
& \frac{2}{5}\!:\!\frac{1}{15}& \frac{1}{40}\!:\!\frac{21}{40}
& \frac{1}{8}& 0\!:\!\tth  
& \frac{1}{8}\!:\!\frac{13}{8}& \tth \!:\!3  \cr
\hline 
\tth & \frac{1}{15}& \frac{21}{40}\!:\!\frac{1}{40}& \frac{1}{15}
& \frac{1}{40}\!:\!\frac{21}{40}
& \frac{2}{5}\!:\!\frac{1}{15}\!:\!\frac{7}{5}& \tth  
& \frac{1}{8}\!:\!\frac{13}{8}& 0\!:\!\tth \!:\!3
& \frac{1}{8}\!:\!\frac{13}{8}  \cr
\hline 
\frac{1}{8}& \frac{1}{40}& \frac{1}{15}\!:\!\frac{2}{5}& \frac{21}{40}
& \frac{1}{15}\!:\!\frac{7}{5}& \frac{1}{40}\!:\!\frac{21}{40}
& \frac{13}{8}& \tth \!:\!3
& \frac{1}{8}\!:\!\frac{13}{8}& 0\!:\!\tth  \cr
\hline 
\end{array} $$
\normalsize

\begin{center}
{Fusion rules of $L(\ff,0)\oplus L(\ff,3)$-modules. }
\end{center}
\small
$$ \begin{array}{|c||c|c|c|c|c|}
\hline
W(0) \! &\!  W(\frac{2}{5}) \! &\!  W(\tth ,+) \! &\! W(\frac{1}{15},+) 
\! &\! 
W(\tth ,-) \! &\! W(\frac{1}{15}.-) \cr
\hline 
\hline
W(\frac{2}{5})\! &\! W(0)\!:\! W(\frac{2}{5})\! &\! 
W(\frac{1}{15},+)\! &\!  W(\frac{1}{15},+)\!:\! 
W(\tth ,+) \! &\!  W(\frac{1}{15},-)\! &\!  W(\frac{1}{15},-)\!:\! 
W(\tth ,-) \cr
\hline
W(\tth ,+) \! &\! W(\frac{1}{15},+)\! &\! W(\tth ,-)\! &\! 
W(\frac{1}{15},-)\! &\! W(0)\! &\! W(\frac{2}{5})   \cr
\hline
W(\frac{1}{15},+) \! &\!   W(\frac{1}{15},+)\!:\! 
W(\tth ,+)\! &\!  W(\frac{1}{15},-)\! &\!  W(\frac{1}{15},-)\!:\! 
W(\tth ,-) \! &\!  W(\frac{2}{5}) \! &\!  W(0)\!:\! W(\frac{2}{5})\cr
\hline
W(\tth ,-)\! &\! W(\frac{1}{15},-)\! &\!  W(0)\! &\!  
W(\frac{2}{5}) \! &\! W(\tth ,+) \! &\! W(\frac{1}{15},+) \cr
\hline
W(\frac{1}{15},-)\! &\!   W(\frac{1}{15},-)\!:\! 
W(\tth ,-)\! &\! W(\frac{2}{5})\! &\!  W(0)\!:\! W(\frac{2}{5}) \! &\! 
W(\frac{1}{15},+)\! &\!  W(\frac{1}{15},+)\!:\! 
W(\tth ,+) \cr
\hline
\end{array} $$

\end{document}